\documentclass[11pt,epsf]{article}
\usepackage{graphicx}
\usepackage{amsthm}
\usepackage{amsfonts}
\usepackage{amssymb}
\title{No zeros of the partial theta function in the unit disk}
\author{Vladimir Petrov Kostov\\
  Universit\'e de Nice, 
Laboratoire de Math\'ematiques, Parc Valrose,\\ 06108 Nice Cedex 2, France,  
e-mail: vladimir.kostov@unice.fr} 
\date{}
\newtheorem{tm}{Theorem}
\newtheorem{defi}[tm]{Definition}
\newtheorem{rem}[tm]{Remark}

\newtheorem{cor}[tm]{Corollary}

\newtheorem{conj}[tm]{Conjecture}
\begin{document}
\maketitle 
\begin{abstract}
  We prove that for $q\in (-1,0)\cup (0,1)$, the partial theta function
  $\theta (q,x):=\sum _{j=0}^{\infty}q^{j(j+1)/2}x^j$ has no zeros in the closed
  unit disk.

{\bf Key words:} partial theta function, Jacobi theta function, 
Jacobi triple product\\ 

{\bf AMS classification:} 26A06
\end{abstract}

\section{Introduction}

We consider the partial theta function
$\theta (q,x):=\sum _{j=0}^{\infty}q^{j(j+1)/2}x^j$, where $q\in [-1,1]$ is a
parameter and $x\in \mathbb{R}$ a variable. In particular,
$\theta (0,x)\equiv 1$, 
$\theta (1,x)=\sum _{j=0}^{\infty}x^j=1/(1-x)$~and

$$\theta (-1,x)=\sum _{j=0}^{\infty}(-1)^jx^{2j}+\sum _{j=0}^{\infty}(-1)^{j+1}x^{2j+1}=
(1-x)/(1+x^2)~.$$ 
For each $q\in (-1,0)\cup (0,1)$ fixed, $\theta$ is an entire function
of $x$ of order~$0$.

The name ``partial theta function'' is connected with the fact that the
{\em Jacobi theta function} equals
$\Theta (q,x):=\sum _{j=-\infty}^{\infty}q^{j^2}x^j$ while
$\theta (q^2,x/q)=\sum _{j=0}^{\infty}q^{j^2}x^j$.  ``Partial'' refers to the fact
that summation in $\theta$ is performed only
  from $0$ to $\infty$. One can observe that 

$$
  \Theta ^*(q,x):=
  \Theta (\sqrt{q},\sqrt{q}x)=\sum _{j=-\infty}^{\infty}q^{j(j+1)/2}x^j=
  \theta (q,x)+\theta (q,1/x)/x~.$$
The function $\theta$ satisfies the relation

  \begin{equation}\label{eqrelation}
    \theta (q,x)=1+qx\theta (q,qx)~.
    \end{equation}
Applications of $\theta$  
to questions concerning asymptotics and modularity of partial and false theta 
functions and their relationship to representation theory and conformal field 
theory (see \cite{CMW} and \cite{BFM}) explain part of the most recent
interest in it. Previously, this function 
has been studied with regard to Ramanujan-type $q$-series 
(see \cite{Wa}), statistical physics 
and combinatorics (see \cite{So}), the theory 
of (mock) modular forms (see \cite{BrFoRh}) and asymptotic analysis 
(see \cite{BeKi}); see also~\cite{AnBe}.

Another domain in which $\theta$ plays an important role is the study of 
section-hyperbolic polynomials. These are real polynomials with
all roots real negative and all whose finite sections (i.e. truncations) 
have also this property, see \cite{KoSh}, \cite{KaLoVi} and \cite{Ost}; 
the cited papers are motivated by results of Hardy, Petrovitch and Hutchinson 
(see \cite{Ha}, \cite{Pe} and \cite{Hu}). Various analytic properties of the
partial theta function are proved in \cite{Ko2}-\cite{Ko12} and other papers
of the author. 

The basic result of the present text is the following theorem (proved in
Section~\ref{secprtmmain}):

\begin{tm}\label{tmmain}
  For each $q\in (-1,0)\cup (0,1)$ fixed, the function $\theta$ has no
  zeros in the closed unit disk~$\overline{\mathbb{D}_1}$.
\end{tm}

In the next section we discuss the question to what extent Theorem~\ref{tmmain}
proposes an optimal result. In Section~\ref{secCOQ}
we make comments and formulate
some open questions.

\section{Optimality of the result}

\subsection{The theorem of Enestr\"om-Kakeya}

For $q\in (0,1)$, the theorem of Enestr\"om-Kakeya about polynomials
with positive coefficients (see \cite{AnSV})
implies that the modulus of each root of
a polynomial $a_0+a_1x+\cdots +a_nx^n$, $a_j>0$, 
is not less than $\min _j|a_{j-1}/a_j|$. When this polynomial equals
$1+qx+\cdots +q^{n(n-1)/2}x^{n-1}$, the minimum equals $1/q$.  
Thus all zeros of all finite truncations of $\theta (q,.)$ (and hence
all zeros of $\theta$ itself) lie outside the
open disk $\mathbb{D}_{1/q}$.

Hence for $q\in (0,1)$
(but not for $q\in (-1,0)$), Theorem~\ref{tmmain}
follows from the theorem of Enestr\"om-Kakeya. 
A hint how to obtain for $q\in (-1,0)$ a disk of a radius tending to $\infty$
as $q\rightarrow 0^-$ and free from zeros of $\theta$
is given in Remark~\ref{remimprove}.

\begin{rem}
  {\rm For $q\in (-1,0)$, it is not true that $\theta (q,.)$ has no zeros
    inside the disk $\mathbb{D}_{1/|q|}$. Indeed, the function $\theta (-0.4,.)$
    has a zero $1.96\ldots <1/0.4=2.5$. More generally, the zero of
    $\theta (q,.)$ closest to the origin can be expanded in a Laurent series
    (convergent for $0<|q|$ sufficiently small) of the form $-1/q-1+O(q)$,
    see~\cite{Ko2}.
    For $q\in (-1,0)$ and $|q|$ sufficiently small, this
    number belongs to the interval $(0,1/|q|)$. See also~\cite{Ko12},
    where the zero set of $\theta$ is
    illustrated by pictures.}
  \end{rem}

\subsection{Optimality w.r.t. the parameter $q$}

(1) This result cannot be generalized in the case when $q$ and $x$ are
    complex. Indeed, suppose that $q\in \mathbb{D}_1$ and $x\in \mathbb{C}$.
    Then the function $\theta$ has no zeros $x$ with $|x|<1/2$. 
    In fact, it has no zero for $|x|\leq 1/2|q|$, see Proposition~7
    in \cite{Ko1}. On the other hand, the radius of the disk in the
    $x$-space centered at $0$ 
in which $\theta$ has no zeros for any $q\in \mathbb{D}_1$ is not larger than 
$0.56\ldots$. Indeed, consider the series $\theta$ with 
$q=\omega :=e^{3\pi i/4}$. 
It equals 

$$(\sum _{j=0}^7\omega ^{j(j+1)/2}x^j)/(1-x^8)~.$$
Its numerator has a simple zero 
$x_*:=0.33\ldots +0.44\ldots i$ whose modulus equals 
$0.56\ldots$. Hence for $\rho \in (0,1)$ sufficiently close to $1$, 
the function $\theta (\rho e^{3i\pi /4},.)$ has a zero close to $x_*$. To see
this one can fix a closed disk $\mathcal{D}$ about $x_*$ of radius $<0.1$
in which $x_*$ is the only zero of $\theta (e^{3i\pi /4},.)$. As $\rho$
tends to $1^-$, the modulus of the difference 
$\theta (e^{3i\pi /4},x)-\theta (\rho e^{3i\pi /4},x)$ tends uniformly to $0$ for
$x\in \partial \mathcal{D}$ (the border of $\mathcal{D}$),
because the series $\theta$ converges uniformly
for $|x|<|x_*|+0.1$, $|q|\leq 1$. The Rouch\'e theorem implies that the
function $\theta (\rho e^{3i\pi /4},.)$ has the same number of zeros in
$\mathcal{D}$ (counted
with multiplicity) as the function $\theta (e^{3i\pi /4},.)$.
\vspace{1mm}

(2) Set $q:=|q|e^{i\phi}$. We show that there exists no interval $J$
on the unit circle centered at $1$ or $-1$ and such that for
$\phi \in J$ and $|q|<1$, the zeros of $\theta (q,.)$
are all of modulus $\geq 1$.  
Suppose that $n\in \mathbb{N}$, $n>2$,
and that $\omega$ is a primitive
root of unity of order $n$. If $n$ is odd, then the sequence of numbers
$\omega ^{k(k+1)/2}$ is $n$-periodic, because $(n+1)/2\in \mathbb{N}$, and
one obtains

$$\theta (\omega ,x)=P(x)/(1-x^n)~,~~~\,
P:=\sum _{j=0}^{n-1}a_jx^j~,~~~\, a_j=\omega ^{j(j+1)/2}~.$$
If $n$ is
even, then this sequence is clearly $(2n)$-periodic, but it is not
$n$-periodic, because $\omega ^{n(n+1)/2}=-1$. One has

$$\theta (\omega ,x)=Q(x)/(1-x^{2n})~,~~~\,
Q:=\sum _{j=0}^{2n-1}b_jx^j~,~~~\, b_j=\omega ^{j(j+1)/2}~.$$
The polynomials $P$ and $Q$ are self-reciprocal, i.~e.
  $a_{(n-1)/2-s}=a_{(n-1)/2+s}$ and
  $b_{(2n-1)/2-s}=b_{(2n-1)/2+s}$. Indeed, for the polynomial $P$ this follows from

$$\begin{array}{rcl}((n-1)/2-s)((n-1)/2-s+1)/2&\equiv&
(n-(n-1)/2+s)(n-(n-1)/2+s-1)/2\\ \\ &=&
((n-1)/2+s)((n-1)/2+s+1)/2\, {\rm mod}[n]~.\end{array}$$
For the polynomial $Q$ one gets

$$\begin{array}{rcl}((2n-1)/2-s)((2n-1)/2-s+1)/2&\equiv&
(2n-(2n-1)/2+s)(2n-(2n-1)/2+s-1)/2\\ \\ &=&
  ((2n-1)/2+s)((2n-1)/2+s+1)/2\, {\rm mod}[2n]~.\end{array}$$
We show that at least one root of the polynomial $P$ and at least one
  root of $Q$ belong to the interior of the unit disk. Indeed, these
  polynomials are monic and $P(0)=Q(0)=1$. The product of their roots being
  equal to $\pm 1$, the only possibility for $P$ and $Q$ not to have roots in
  $\mathbb{D}_1$ is all their roots to be of modulus~$1$. These polynomials are
  self-reciprocal, so $P(z)=0$ (resp. $Q(z)=0$) implies $P(1/z)=0$
  (resp. $Q(1/z)=0$). But if $|z|=1$, then $1/z=\bar{z}$. This means that $P$
  and $Q$ can have as roots either $\pm 1$ or complex conjugate pairs, i.~e.
  $P$ and $Q$ must be real which is false as their coefficients of $x$ equal
  $\omega \neq \pm 1$.

  So $P$ and $Q$ have each at least one root in $\mathbb{D}_1$. As in part (1)
  of this subsection one deduces that for $|q|$ sufficiently close to $1$ and
  for $e^{i\phi}=\omega$, the function $\theta (q,.)$ has a zero in
  $\mathbb{D}_1$. Primitive roots are everywhere dense on the unit circle.
  This implies the
  absence of an interval~$J$ as above.

  \subsection{Optimality w.r.t. the variable $x$}

  Suppose first that $q\in (-1,0)$. Then in the formulation of
  Theorem~\ref{tmmain} one cannot replace the unit disk
by a disk of larger radius. Indeed, the zero
of the numerator of $\theta (-1,x)$ (which equals $1$) is the limit
as $q$ tends to $-1^+$ of the smallest positive zero of $\theta (q,x)$, see
\cite[Part~(2) of Theorem~3]{Ko12}, so in any disk $\mathbb{D}_{1+\varepsilon}$,
$\varepsilon >0$, there is a zero of $\theta$ for some $q\in (-1,0)$.

Suppose now that $q\in (0,1)$.

\begin{conj}
  Theorem~\ref{tmmain} does not hold true if one replaces in its formulation
  the unit disk by a disk of larger radius.
  \end{conj}

The following numerical example shows why this conjecture should be
considered as plausible. Set $\theta _{100}:=\sum _{j=0}^{100}q^{j(j+1)/2}x^j$
(the $100$th truncation of $\theta$).
For $q=0.98$, the function $\theta _{100}(0.98,.)$ has a zero
$\lambda _0:=1.209\ldots +0.511\ldots i$, of modulus
$1.312\ldots$. For $q=0.98$ and
$x=1.32$, the first two terms of $\theta$ which are not in $\theta _{100}$
equal $y_{101}:=7.407\ldots \times 10^{-33}$ and
$y_{102}:=1.270\ldots \times 10^{-33}$
respectively. Their ratio is $y_{101}/y_{102}>5.5$ and the moduli of the terms of
$\theta$ decrease faster than a geometric progression.
Hence for $|x|<1.32$, one has 

$$T_0:=|\theta (0.98,x)-\theta _{100}(0.98,x)|<
y_{101}/(1-5.5^{-1})=9.053\ldots \times 10^{-33}~.$$
On the other hand
$\Lambda _0:=(\partial \theta /\partial x)(0.98,\lambda _0)=27.180\ldots +
18.959\ldots i$, with $|\Lambda _0|>33$. Thus one should expect to
find a zero of $\theta (0.98,.)$ close to $\lambda _0$ (the truncated terms
are expected to change the position of $\lambda _0$ by
$\approx T_0/|\Lambda _0|$ which quantity is of order $10^{-34}$.
So in the formulation of Theorem~\ref{tmmain} one should not be able to replace
the unit disk by a disk of radius larger than~$1.32$.

\section{Proof of Theorem~\protect\ref{tmmain}\protect\label{secprtmmain}}

We remind first that the {\em Jacobi triple product} is the identity

$$\Theta (q,x^2)=\prod_{m=1}^{\infty}(1-q^{2m})
       (1+x^2q^{2m-1})(1+x^{-2}q^{2m-1})$$
       which implies
       $\Theta ^*(q,x)=\prod_{m=1}^{\infty}(1-q^m)(1+xq^m)(1+q^{m-1}/x)$. 
       Thus

       \begin{equation}\label{eqT}
         \prod_{m=1}^{\infty}(1-q^m)(1+xq^m)(1+q^{m-1}/x)=
         \theta (q,x)+\theta (q,1/x)/x~.
         \end{equation}

  Suppose that $q\in (-1,0)\cup (0,1)$, $x_0\in \mathbb{C}$, $|x_0|=1$
  (hence $\overline{x_0}=1/x_0$), and
  that $\theta (q,x_0)=0$. The coefficients of $\theta$ being real, one has
  $\theta (q,\overline{x_0})=\overline{\theta (q,x_0)}=0$, so the right-hand
  side of equation (\ref{eqT}) equals~$0$ for $x=x_0$. However for $x=x_0$, the
  left-hand side vanishes only for $x_0=-1$.

  For $q\in (0,1)$, one has $\theta (q,-1)=\sum _{j=0}^{\infty}(-1)^jq^{j(j+1)/2}$,
  and the
  latter function takes only values from the interval $(1/2,1)$, with
  $\lim _{q\rightarrow 1^-}=1/2$, see~\cite[Propositions~14 and~16]{Ko1}. For
  $q\in (-1,0)$, one sets $u:=-q$, so

  $$\begin{array}{rcl}
    \theta (q,-1)=\theta (-u,-1)&=&1+u-u^3-u^6+u^{10}+u^{15}-u^{21}-u^{28}+
    \cdots \\ \\ &=&
  1-u^3+u^{10}-u^{21}+\cdots +u-u^6+u^{15}-u^{28}+\cdots >0~,\end{array}$$
  because this is the sum of two Leibniz series with positive initial terms.
  Thus for $q\in (-1,0)\cup (0,1)$, the partial theta function has no zeros
  of modulus~$1$.

  For $-1/2<q<1/2$, one has $\theta (q,x)\neq 0$ for any
  $x\in \overline{\mathbb{D}_1}$, because

  $$|\theta (q,x)|\geq 1-|q|-|q|^3-|q|^6-\cdots \geq 
  1-|q|-|q|^2-|q|^3-\cdots =(1-2|q|)/(1-|q|)>0~.$$
  As the parameter $q$ varies in $(0,1)$ or in $(-1,0)$, the zeros of $\theta$
  depend continuously on $q$. For $|q|<1/2$, there are no zeros of $\theta$
  in $\overline{\mathbb{D}_1}$ and for $q\in (-1,0)\cup (0,1)$, no zero of
  $\theta$ crosses $\partial \mathbb{D}_1$
  (the border of the unit disk). Hence for
  $q\in (-1,0)\cup (0,1)$, there are no zeros
  of $\theta$ in $\overline{\mathbb{D}_1}$.

  \begin{rem}\label{remimprove}
    {\rm One can prove that for $|q|\leq 0.4$, the function $\theta (q,.)$
      has no zeros in the closed disk $\overline{\mathbb{D}_{1/\sqrt{|q|}}}$.
      Indeed,}
    $$|\theta (q,1/\sqrt{|q|})|\geq 1-\sum _{j=1}^{\infty}|q|^{j(j+1)/2-j/2}=
    1-\sum _{j=1}^{\infty}|q|^{j^2/2}\geq 
    1-\sum _{j=1}^{\infty}0.4^{j^2/2}=0.19\ldots >0~.$$
    \end{rem}

  \section{Comments and open questions\protect\label{secCOQ}}

  \subsection{The case $q\in (0,1)$.}

\begin{figure}[htbp]
%\includegraphics[scale=0.5]{parthetanegfirstfour.eps}
%\centerline{\hbox{\includegraphics[scale=0.7]{parthetanegfirstfour.eps}}}
\centerline{\hbox{\includegraphics[scale=0.7]{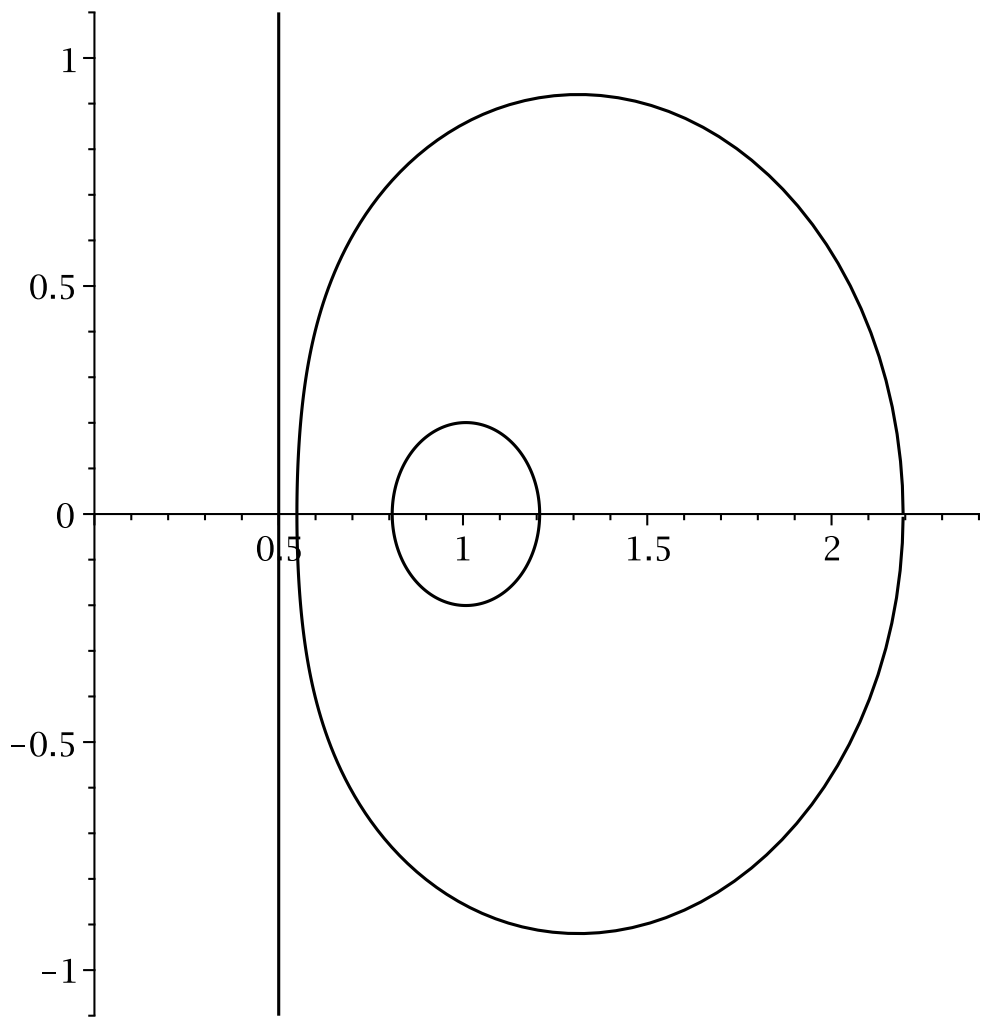}}}
%\centerline{\hbox{\epsfxsize=10cm \epsfbox{k=1234.pdf}}}
\caption{The vertical line Re$x=1/2$ and the images of the unit circle
  under the map $x\mapsto \theta (q,x)$
for $q=0.2$ and $q=0.7$.}
\label{imageunitdiskqpos}
\end{figure}

  In Fig.~\ref{imageunitdiskqpos} we show the images for $q=0.2$
  (the smaller oval) and for $q=0.7$ (the larger oval) 
  of the unit circle in the $x$-plane
  under the mapping $x\mapsto \theta (q,x)$, together with the vertical line
  Re$x=1/2$.  

It would be interesting to know whether:
\vspace{1mm}

1) for any $q\in (0,1)$, the image of the unit circle is a convex oval
about the point $(1,0)$ and 
belonging to the half-plane Re$x>1/2$;
\vspace{1mm}

2) for $0<q_1<q_2<1$, the image of the unit circle for $q=q_1$ lies
inside its image for $q=q_2$.
\vspace{1mm}

These questions are motivated by the fact that for $q=1$, one
has $\theta (1,x)=1/(1-x)$, and for $|x|=1$, it is true that
Re$(1/(1-x))=1/2$, i.e. the vertical line Re$(1/(1-x))=1/2$
is the image of the unit circle for $q=1$; on the other hand, the point $(1,0)$
is the image of the unit circle for $q=0$.

\subsection{The case $q\in (-1,0)$.}

\begin{figure}[htbp]
%\includegraphics[scale=0.5]{parthetanegfirstfour.eps}
%\centerline{\hbox{\includegraphics[scale=0.7]{parthetanegfirstfour.eps}}}
\centerline{\hbox{\includegraphics[scale=0.7]{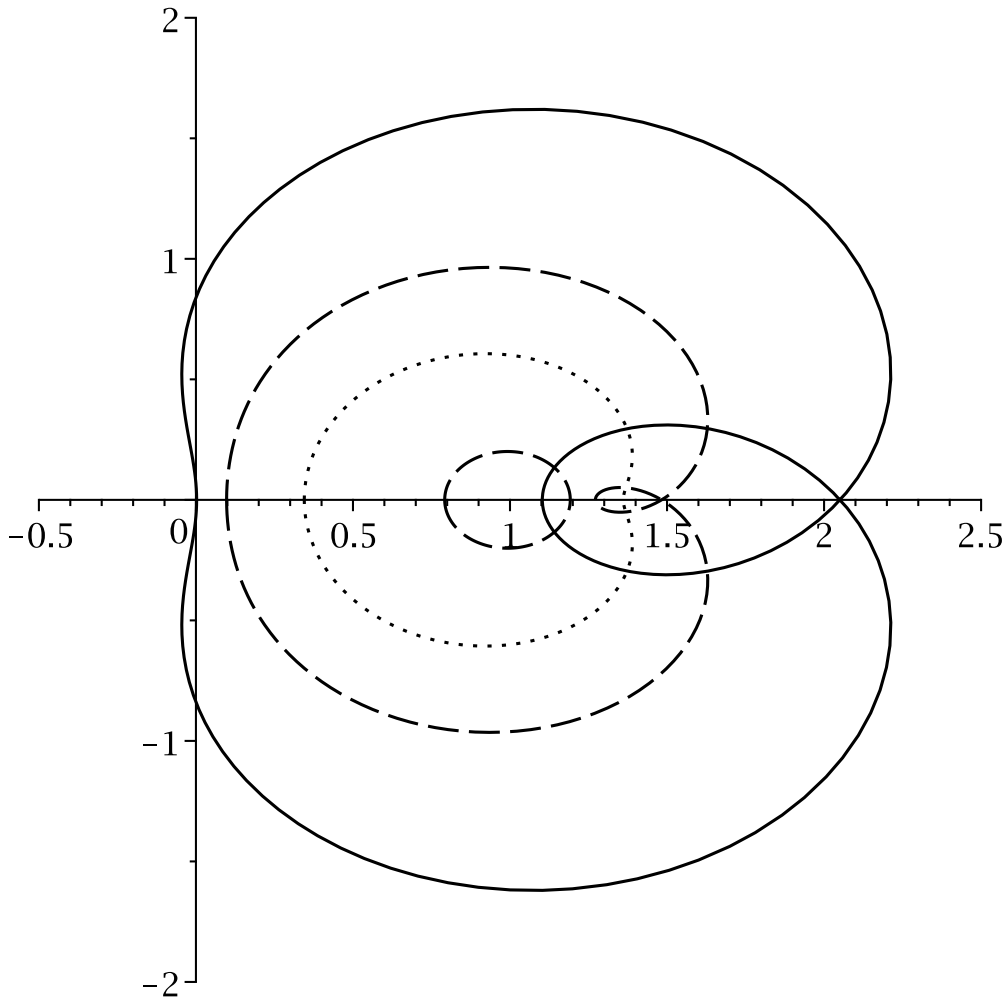}}}
%\centerline{\hbox{\epsfxsize=10cm \epsfbox{k=1234.pdf}}}
\caption{The images of the unit circle
  under the map $x\mapsto \theta (q,x)$
for $q=-0.2$, $q=-0.53$, $q=-0.7$ and $q=-0.85$.}
\label{imageunitdiskqneg}
\end{figure}

In Fig.~\ref{imageunitdiskqneg} we show the images for $q=-0.2$
(small oval in dashed line), $q=-0.53$ (closed contour in dotted line),
$q=-0.7$ (curve with self-intersection in dashed line) and $q=-0.85$
(curve with self-intersection in solid line) 
  of the unit circle in the $x$-plane
  under the mapping $x\mapsto \theta (q,x)$.
  
\begin{figure}[htbp]
%\includegraphics[scale=0.5]{parthetanegfirstfour.eps}
%\centerline{\hbox{\includegraphics[scale=0.7]{parthetanegfirstfour.eps}}}
\centerline{\hbox{\includegraphics[scale=0.7]{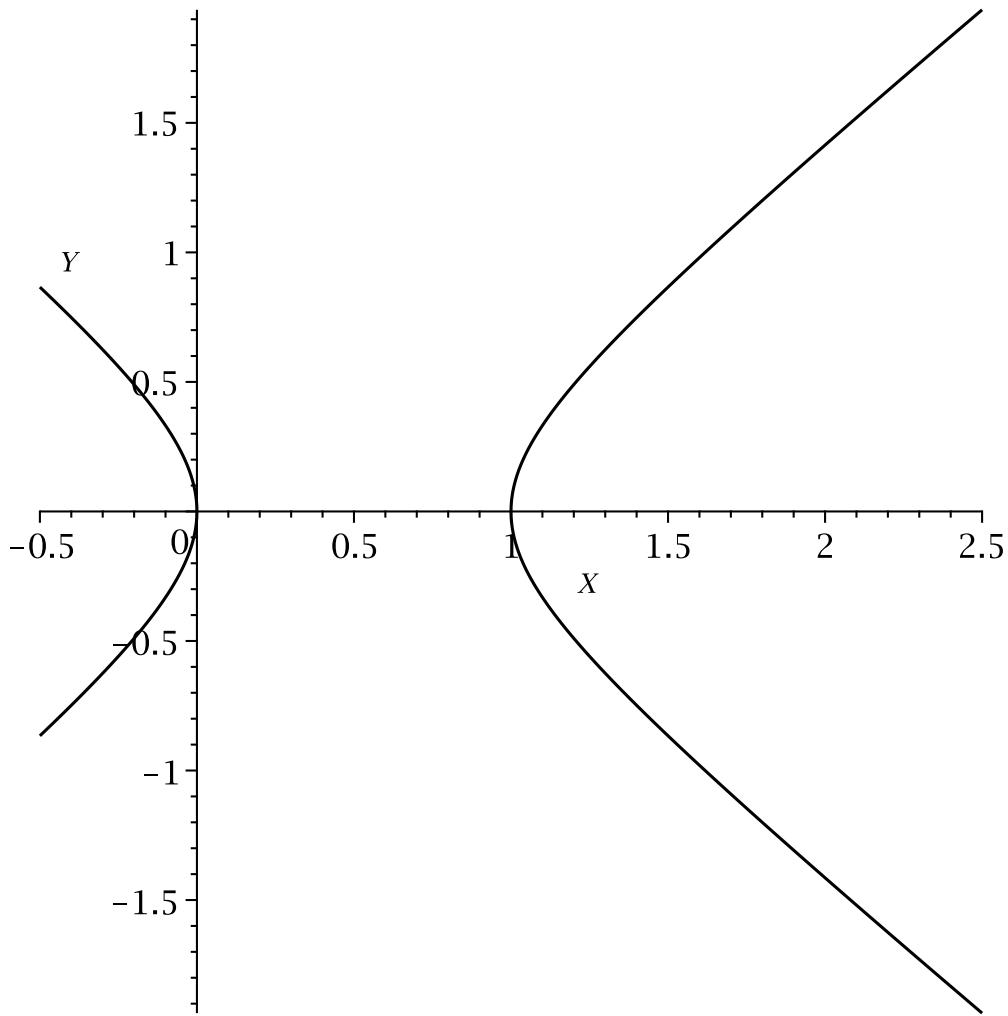}}}
%\centerline{\hbox{\epsfxsize=10cm \epsfbox{k=1234.pdf}}}
\caption{The image of the unit circle
  under the map $x\mapsto \theta (q,x)$
for $q=-1$.}
\label{imageunitdiskqnegbis}
\end{figure}

The following questions are natural to ask:
\vspace{1mm}

3) Is it true, and for which value of $v\in (0,1)$, that for $q\in (-v,0)$,
the corresponding image is a convex oval about the point $(1,0)$?
\vspace{1mm}

4) Is it true that for $q\in (-w,-v)$, $-1<-w<-v<0$, the corresponding image
changes convexity twice ``at its right'' (as
this seems to be the case of the curve given in dotted line)?
\vspace{1mm}

5) Is it true that for $q\in (-1,-w)$, the image has a self-intersection
point? One can expect that for $q=-w$, the image has a cusp point.
\vspace{1mm}

6) Is it true that for $q\in (-1,-w')$, $-1<-w'<-w$, the image has
still self-intersection and changes convexity twice ``at its left''?
\vspace{1mm}

7) Is it true that for $q\in (-1,-w'')$, $-1<-w''<-w'$, the image has
still self-intersection, changes convexity twice ``at its left''
and intersects the vertical axis at four points?
(For $q=w''$, the image is supposed to have two
tangencies with the vertical axis.)
\vspace{1mm}

8) Is it true that these are all transformations which the image undergoes
for $q\in (-1,0)$?
\vspace{1mm}

9) Is it true that for $-1<q_2<q_1<0$, the image of the unit circle for
$q=q_1$ lies inside its image for $q=q_2$? ``Inside'' means ``inside the
contour excluding (for $q_2>w$) the loop''.
\vspace{1mm}

It should be observed that for values of $q$ close to $-1$, the image seems to
pass through the origin. In reality, it passes very close to it, but
nevertheless to its right, according to Theorem~\ref{tmmain}.
The image of the unit circle for $q=-1$ is the hyperbola $Y^2-X^2-X=0$, where
$X:=$Re$x$ and $Y:=$Im$x$, see Fig.~\ref{imageunitdiskqnegbis}. (The centre
of the hyperbola is at $(1/2,0)$, its asymptotes are the lines
$Y=\pm (X-1/2)$.) Following
a similar logic one can assume that the point $(1,0)$ remains in the exterior
of the loop of the image (the loop existing for $q>w$). The proximity
of the image to
the origin makes it seem unlikely that one could prove the absence of zeros
of $\theta$ in a disk of a radius larger than $1$ (for all $q\in (-1,0)$). 
\vspace{1mm}

{\bf Acknowledgement.} The author is grateful to A.~Vishnyakova and B.~Shapiro 
for the useful comments of this text.

\end{document}